# Continuous data assimilation for problems with limited regularity using non-interpolant observables

Vladimir Yushutin[*]

**Abstract**

Continuous data assimilation addresses time-dependent problems with unknown initial conditions by incorporating observations of the solution into a nudging term. For the prototypical heat equation with variable conductivity and the Neumann boundary condition, we consider data assimilation schemes with non-interpolant observables unlike previous studies. These generalized nudging strategies are notably useful for problems which possess limited or even no additional regularity beyond the minimal framework. We demonstrate that a spatially discretized nudged solution converges exponentially fast in time to the true solution with the rate guaranteed by the choice of the nudging strategy independent of the discretization. Furthermore, the long-term discrete error is optimal as it matches the estimates available for problems of limited regularity with known initial conditions. Three particular strategies – nudging by a conforming finite element subspace, nudging by piecewise constants on the boundary mesh, and nudging by the mean value – are explored numerically for three test cases, including a problem with Dirac delta forcing and the Kellogg problem with discontinuous conductivity.

## 1 Introduction

Consider the following weak formulation of the heat equation on $(0, T) \times \Omega$ for a bounded Lipschitz domain $\Omega \subset \mathbb{R}^n$ with the inhomogeneous Neumann boundary condition on $\partial \Omega$: given $g \in L^2(0, T; H^{-1/2}(\partial \Omega))$ and $f \in L^2(0, T; H^{-1}(\Omega))$, find $\nu \in L^2(0, T; H^1(\Omega)) \cap H^1(0, T; H^{-1}(\Omega))$ such that

$$\langle \partial_t \nu, v \rangle + a(\nu, v) = \langle f, v \rangle + \langle g, v \rangle_\partial, \qquad \forall v \in H^1(\Omega), \tag{1.1}$$

where a uniformly elliptic bilinear form $a(u, v)$ is given by $(A \nabla u, \nabla v)$ with a symmetric, uniformly Lipschitz matrix $A$. Here $(\cdot, \cdot)$ denotes the $L^2(\Omega)$ inner product, and the notations for other relevant Sobolev norms and products are introduced as follows,

$$\begin{aligned}
\|u\|_0^2 = \|u\|_{L^2(\Omega)}^2, \quad |u|_a^2 = a(u, u), \quad \|u\|_1^2 = \|u\|_0^2 + |u|_a^2, \quad \|u\|_\partial^2 = \|u\|_{L^2(\partial\Omega)}^2, \\
\langle g, v \rangle_\partial = \langle g, v \rangle_{H^{-1/2}(\partial\Omega), H^{1/2}(\partial\Omega)}, \quad \langle f, v \rangle = \langle f, v \rangle_{H^{-1}(\Omega), H^1(\Omega)}.
\end{aligned} \tag{1.2}$$

The problem (1.1) needs an initial condition for some $\nu_0 \in L^2(\Omega)$,

$$\nu(0) = \nu_0. \tag{1.3}$$

---

[*]Department of Mathematics, University of Tennessee, Knoxville. Email: vyushuti@utk.edu

However, $\nu_0$ is not available in some situations. For instance, in modeling of the Earth's mantle the initial state is inherently illusive, see [Van+22, Section 5.3.2]. To circumvent the issue, one may opt to conduct *continuous data assimilation* [AOT14] by replacing the unknown initial condition $\nu_0$ with a convenient $u_0 = 0$ and composing the following nudged variational problem,

$$\langle \partial_t u, v \rangle + a(u,v) + \mu(L_H u - \nu_H, L_H v)_Y = \langle f, v \rangle + (g,v)_\partial, \quad \forall v \in H^1(\Omega), \quad u(0) = 0, \quad (1.4)$$

whose solution $u(t)$ is expected to converge to the exact solution $\nu(t)$ as $t$ evolves thanks to a *nudging* term with the parameter $\mu$. The missing initial data $\nu_0$ in (1.3) is compensated by continuous in time *observations* $\nu_H(t)$ in (1.4) during an initial stage $t \in (0, T_\mu)$, for some $T_\mu \leq T$. The observations $\nu_H(t)$ of exact solution $\nu(t)$ are obtained using an *observation mapping* $L_H$,

$$\nu_H = L_H \nu, \qquad L_H : H^1(\Omega) \to V_H, \quad (1.5)$$

where $V_H$ is a Hilbert space equipped with $(\cdot, \cdot)_Y$ inner product. The basic requirement for the linear operator $L_H$ in (1.5) is its continuity, i.e. there exists a constant $C_B > 0$ such that

$$C_B^{-1} \|L_H u\|_Y^2 \leq \|u\|_{L^2(\Omega)}^2 + a(u,u), \qquad \forall u \in H^1(\Omega). \quad (1.6)$$

Assume now that we introduce a conforming, spatial discretization by $V_h \subset H^1(\Omega)$ of the problem (1.4). From the computational point of view, the convergence of the nudged solution $u_h(t)$ of a spatially discretized scheme for (1.4) to the exact $\nu(t)$ is preferred to be exponentially fast in time,

$$\|u_h(t) - \nu(t)\|_0 \lesssim e^{-t/T_\mu}, \quad (1.7)$$

with the characteristic time $T_\mu$ such that $T_\mu \ll T$. In this case, $u_h(T_\mu)$ can serve as the restarting initial condition, $\nu_h(T_\mu) = u_h(T_\mu)$, in the discretization of the heat equation (1.1) on $(T_\mu, T)$,

$$\langle \partial_t \nu_h, v_h \rangle + a(\nu_h, v_h) = \langle f, v_h \rangle + \langle g, v \rangle_\partial, \quad \forall v_h \in V_h, \quad (1.8)$$

provided the error $\|\nu(T_\mu) - u_h(T_\mu)\|_0$ of nudging falls below the unavoidable error in the spatial approximation of the initial condition $\nu(T_\mu)$ by $V_h$. Ideally, the initial stage $(0, T_\mu)$ in (1.7) shortens as the nudging parameter $\mu$ in (1.4) increases. We will see that actually this is not the case generally, and the practical choice of $\mu$ is limited from above by the so-called *saturated* value (1.12) of the nudging parameter.

## 1.1. Motivation and outline

The idea of restarting (1.8) at $t = T_\mu$ requires the convergence of the numerical solution $u_h(t)$ in the $L^\infty(0, T_\mu; L^2(\Omega))$ norm which is not expected even for the discrete solution $\nu_h$ of (1.8) with a proper initial condition $\nu_h(0)$ if the problem (1.1) is not sufficiently regular. The first goal of this paper is to reconcile this situation and develop a nudging approach suitable for minimally regular problems. To be specific, we will prove that if the spatially discretized nudged problem (1.4) is solved on the whole $(0, T)$, then, after the period $(0, T_\mu)$ of exponential decay (1.7), its solution $u_h$ converges to the exact $\nu$ in the $L^2(T_\mu, T; L^2(\Omega))$ sense. The developed theory suggests the optimal rates of the error of nudging that match what is expected from the non-nudged scheme (1.8) with proper initial condition in case of a problem with limited regularity, i.e. if $\nu \in L^2(0, T, H^{1+s}(\Omega))$, then

$$\|u_h - \nu\|_{L^2(T_\mu, T; L^2(\Omega))} \lesssim h^{s+1}. \quad (1.9)$$

In addition, we consider the heat equation (1.1) with the variable conductivity which appears to be new in the context of continuous data assimilation.

The second motivation stems from the fact that not all observations that may be useful for data assimilation involve an approximation subspace of $H^1(\Omega)$. Therefore, an approximation error estimate $\|L_H u - u\|_0 \lesssim H\|u\|_1$, $L_H u \in H^1(\Omega)$, may not hold for a general strategy (1.5). For instance, if one simulates the Earth's mantle convection, almost all historical observations available for nudging are sourced from the lithosphere which is rather the boundary of the domain of interest, see in [Ger14, Figure 1]. To broaden the scope of possible ways to incorporate the observational data, we consider the abstract setting in which the observation space $V_H$ is not related to $H^1(\Omega)$ by any means other than the mapping $L_H$ in (1.5). Moreover, the abstract setting becomes even more important in the search of a nudging strategy for minimally regular problems; for example, any point-based operator such as Lagrange interpolation is not well-defined in this case. In this regard, results of this paper hold for any nudging strategy (1.5) provided the abstract Assumtpion 1.1 and Assumption 2.1 are satisfied.

We would like to emphasize that the paper focuses on the spatial discretization only because it is seems to be critical to the outlined goals. The full discretization should be analyzed elsewhere and the results of this paper can be seen as the best case scenario.

The paper is organized as follows. In the remaining part of Section 1 we introduce the observability Assumption 1.1 and connect it with the so-called *saturation effect*. Three examples of nudging strategies are discussed along with correspondent saturated choices of $\mu$ such that the $u(t)$ solving (1.4) converges in time to $\nu(t)$ solving (1.1) as $e^{-\mu t}$ in $L^2(\Omega)$ norm. The latter is proven in Proposition 2.4 of Section 2 which is devoted to the well-posedness and regularity of the continuous nudged problem (1.4). Note that the abstract stability Assumption 2.2 is introduced in Section 2 and it backbones the semi-discrete error analysis of nudging schemes (3.1) presented in Section 3. An extensive numerical investigation of three examples of nudging strategies from Section 1.4 is conducted in Section 4. In addition to a smooth test case, we nudge toward a minimally regular solution of (1.1), and even consider the Kellogg problem with discontinuous coefficients, see Section 4.1. The obtained numerical results verify the optimal error analysis of Section 3 and showcase the saturation effect.

### 1.2. Literature review

A literature review aligned with the focus of this paper is in order. Among the papers on nudging that are related directly to the numerical analysis of prototypical heat problem (1.1), [DLR25] considers the case of constant conductivity and the homogeneous Dirichlet boundary condition, and the authors prove optimal error estimates assuming $\nu \in L^\infty(0,T;H^2(\Omega))$. In [RZ21], the linear advection-diffusion problem with the same assumed regularity is analyzed as well. Surprisingly, the guaranteed exponential decay in both papers is bounded from above by the Poincaré constant and independent of the approximability of $V_H$ while the numerical experiments suggest quite the opposite. Historically, the research on continuous data assimilation has been focused on the Navier-Stokes equation driven by the goal of predicting weather, see [BM17, Section 1]. The discrete error analysis of related schemes is conducted in [LRZ19; GN20; IMT20; GNT20; RZ21; DLR25] only and the common to them regularity assumption is $L^\infty(0,T;H^2(\Omega))$. Another work that includes numerical analysis is [DR22] for $L^\infty(0,T;H^3(\Omega))$ solutions of Cahn-Hilliard problem. Note that there are PDEs for which continuous data assimilation is not applicable at all [TV24] as they do not possess a finite number of so-called determining parameters [FP67; AOT14].

Existing papers on continuous data assimilation seem to rely solely on interpolation operators, i.e. considering $L_H := I_H$ in the nudging term of (1.4), with the approximation property $\|I_H u - u\|_0 \lesssim H \|u\|_1$. Interesting enough, the foundational paper [AOT14] was accompanied by [AT14], Section 6 of which is devoted to a non-approximation, distributional nudging strategy; however, this direction does not seem to have been explored further. We now mention papers that involve observations that are still approximational, but are utilized in non-standard ways. First of all, the $L^2$-projection operator onto the space $V_H$ of piecewise constants is perturbed on the discrete level in [RZ21] by lumping its matrix. This lumped version is simple to implement but only in the case of observational nodes being the nodes of the fine mesh as well. A similar alignment of the coarse and fine meshes leads to the $\mu$-independent error estimates in [DLR25]; also see Corollary 3.3 related to it. In [Vic+21], the authors consider observational nodes that move in space and demonstrate the feasibility of nudging in computations. A time-averaged interpolant operator is studied in [Jol+19] providing no error estimates. The authors of [BBJ21] consider the nudging strategy for a spectral method using observations on a moving subdomain of $\Omega$.

### 1.3. Spatial observability and the saturation effect

We first introduce an assumption on the observational mapping $L_H$ which will play the central role in the continuous analysis of the nudged parabolic problem in Section 2 as well as in the numerical analysis of it in Section 3. For instance, the *spatial observability* property 1.10 will be crucial for Proposition 2.4 that guarantees exponential convergence in time of the nudged solution of (1.4) to the exact one of (1.1). Before we proceed, it is worth to separate the *spatial observability* (1.10) from the *observability properties* well-known in the control theory for PDEs, see e.g. [Zua07, Section 5]. First, these observability inequalities estimate the adjoint state by the $L^2(0,T;L^2(\omega))$ norm of observations, where $\omega \subset \Omega$ is domain involved in the *interior control* problem. Second, the inequalities hold for the solution, not any function as in (1.10).

**Assumption 1.1** (Observability)**.** For the operator $L_H : H^1(\Omega) \to V_H$ and the space $V_H$ with the norm $Y$ in (1.5), there exist a constant $C_a \geq 0$ such that

$$\|u\|_0^2 \leq \|L_H u\|_Y^2 + C_a |u|_a^2, \quad \forall u \in H^1(\Omega). \tag{1.10}$$

The property establishes the control of the state by the observations with the help of the $|u|_a^2$ term, and a smaller constant $C_a$ indicates a better nudging strategy. We will see in Proposition 2.4 that the guaranteed rate in $e^{-\gamma t}$ of the exponential decay in time of the nudged problem with the nudging parameter $\mu$ is given by

$$\gamma = \min\left(\mu, \frac{1}{C_a}\right) \leq \gamma_s := C_a^{-1}, \tag{1.11}$$

where the rate $\gamma_s$ is attained by choosing in (1.4) any value of $\mu$ such that

$$\mu > \mu_s := C_a^{-1} \tag{1.12}$$

For a fixed nudging strategy (1.5), any value of $\mu$ larger than $\mu_s$ (1.12) cannot guarantee the increase of $\gamma$ in Proposition 2.4. For this reason one calls $\mu_s$ and $\gamma_s$ above the *saturated* nudging parameter and saturated rate, correspondingly.

*Remark* 1.2. We would like to emphasize that Assumption 1.1 holds for any $u \in H^1(\Omega)$ with the same constant $C_a$. However, for a specific $u \in H^1(\Omega)$ the inequality (1.10) holds with the constant $C_a[u] = (\|u\|_0^2 - \|L_H u\|_Y^2)/|u|_a^2 < C_a$, $|u|_a \neq 0$, which may be substantially smaller than $C_a$, e.g. for a $u$ with high gradients so $|u|_a$ is large relative to other norms in (1.10). In this sense, nudging to a solution of (1.1) with higher gradients is expected to be faster.

As can be seen from (1.11), the saturation effect depends critically on the observability property (1.10) and, therefore, on the nudging strategy consisting of the choice of mapping $L_H$, space $V_H$ and the norm $Y$ in (1.5). In the next section we present three examples of nudging strategies and discuss practical choices of the nudging parameter in view of the saturation effect.

### 1.4. Examples of nudging strategies

Here we present three straightforward examples of an observational mapping $L_H$, a space $V_H$ and a norm $Y$ satisfying Assumption 1.1 that will be utilized later in numerical experiments of Section 4. Note that throughout this section the symbol $\lesssim$ hides the constants stemming from the domain $\Omega$ and, if $V_H$ involves a triangulation, from its quality. In the rest of the paper these constants are embedded in $C_a$ which will be always mentioned explicitly.

- **Nudging by an approximating projection.** Let the space $V_H := Q_1(\mathcal{T}_H) \subset H^1(\Omega)$ be a conforming finite element subspace of degree 1 based on a quasi-uniform and shape-regular quadrilateral triangulation $\mathcal{T}_H$ of $\Omega$ with the mesh size $H$. The observation mapping $L_H$ is given by the $L^2$-projection, $P_H : H^1(\Omega) \to Q_1(\mathcal{T}_H)$, defined weakly by

$$(P_H u, r) = (u, r), \qquad \forall r \in Q_1(\mathcal{T}_H).$$

The continuity 1.6 is satisfied trivially for the choice $Y := \|\cdot\|_0$ of the nudging norm. Moreover, the observability property (1.10) can be established with $C_a \simeq H^2$,

$$\|u\|_0^2 \leq \|P_H u\|_0^2 + \|P_H u - u\|_0^2 \lesssim \|P_H u\|_0^2 + H^2 |u|_a^2, \quad \forall u \in H^1(\Omega), \qquad (1.13)$$

using the well-known approximation error estimate, $\|P_H u - u\|_0 \lesssim H|u|_a$, of the projection $P_H$ onto the finite element space $V_H(\mathcal{T}_H)$. Consequently, according to Proposition 2.4, the nudging strategy given by

$$V_H := Q_1(\mathcal{T}_H), \qquad L_H := P_H, \qquad Y := \|\cdot\|_0 \qquad (1.14)$$

guarantees the rate (1.11) of the exponential decay of the error of nudging $u - \nu$ given in

$$\gamma_s = C_a^{-1} \simeq H^{-2}, \qquad \mu_s \simeq H^{-2} \qquad (1.15)$$

along with a practical choice for $\mu$ in (1.4).

- **Nudging by the boundary projection.** Consider the finite element space of degree 0, $V_H := Q_0(\partial \mathcal{T}_H) \subset L^2(\partial \Omega)$, consisting of piecewise constants on a quasi-uniform, shape-regular quadrilateral triangulation $\partial \mathcal{T}_H$ of $\partial \Omega$ with the mesh size $H$. The observation operator $L_H$ can be chosen to be the projection, $P_H^\partial : L^2(\partial \Omega) \to Q_0(\partial \mathcal{T}_H)$, defined weakly by

$$(P_H^\partial u, r)_\partial = (u, r)_\partial, \qquad \forall r \in Q_0(\partial \mathcal{T}_H),$$

where $(u, v)_\partial = \int_{\partial \Omega} uv$ and (1.6) is satisfied trivially for the choice $Y := \|\cdot\|_\partial$ of the nudging norm. To show that the property (1.10) holds with $C_a \simeq 1$, we invoke the Poincare-Steklov inequality in the following estimate,

$$\|u\|_0^2 \lesssim \|u\|_\partial^2 + |u|_a^2 \lesssim \|\Pi_H^\partial u\|_\partial^2 + H\|u\|_0^2 + (1+H)|u|_a^2, \quad \forall u \in H^1(\Omega),$$

in which the trace theorem and errors estimates for the boundary projection yielded $\|\Pi_H^\partial u - u\|_\partial^2 \lesssim H\|u\|_1^2$. Consequently, for a sufficiently small $H$ the nudging strategy

$$V_H := Q_0(\partial \mathcal{T}_H), \qquad L_H := P_H^\partial, \qquad Y := \|\cdot\|_\partial \qquad (1.16)$$

guarantees the rate (1.11) of the exponential decay of the error of nudging $u - \nu$ given in

$$\gamma_s \simeq C_a^{-1} \simeq 1, \qquad \mu_s \simeq 1 \tag{1.17}$$

along with a practical choice for $\mu$ in (1.4).

- **Nudging by the mean value.** For parabolic problem with Neumann boundary condition, the evolution in time of the mean value provides valuable information for nudging. In this case, the observability property (1.10) holds with the Poincaré constant, $C_a := C_P$,

$$\|u\|_0^2 \leq \|\bar{u}\|_0^2 + \|u - \bar{u}\|_0^2 \leq |\Omega|^{-1}\bar{u}^2 + C_P|u|_a^2,$$

where $\bar{u} := |\Omega|^{-1} \int_\Omega u$ is the mean value. We see that the nudging strategy

$$V_H := \mathbb{R}, \qquad L_H u := \bar{u}, \qquad Y := |\Omega|^{1/2}|\cdot| \tag{1.18}$$

guarantees the rate (1.11) of the exponential decay of the error of nudging $u - \nu$ given in

$$\gamma_s = C_P^{-1}, \qquad \mu_s = C_P^{-1} \tag{1.19}$$

along with a practical choice for $\mu$ in (1.4).

## 2 Nudged problem

The nudged problem (1.4) is the weak form of the following equation,

$$\partial_t u + (\Delta_a + \mu L_H^* L_H)u = f + \mu L_H^* \nu_H, \qquad \partial_n u = g, \qquad u(0) = 0, \tag{2.1}$$

where $\Delta_a u := -\operatorname{div}(A\nabla u)$ is the variable-coefficient Laplacian with conductivity matrix $A$ and $L_H^* : V_H \to H^{-1}(\Omega)$ is the adjoint to (1.5) operator defined for any $z \in V_H$ by

$$\langle L_H^* z, v \rangle = (z, L_H v)_Y, \qquad \forall v \in H^1(\Omega).$$

Thanks to (1.6), the mapping $L_H^*$ in bounded as well. Note that the term $L_H^* \nu_H$ may be less regular than the forcing $f$, therefore, the nudged solution $u$ may be less regular than the exact solution $\nu$. For this reason, a stronger condition on the nudging strategy is assumed following the existing literature (see Section 1.2) on the nudging via interpolation observables.

**Assumption 2.1** (Stability). The operator $L_H : H^1(\Omega) \to V_H$ is stable in the $L^2(\Omega)$ norm, i.e. there exists a constant $C_S > 0$ which may depend on $V_H$ such that

$$\|L_H u\|_Y^2 \leq C_S \|u\|_0^2, \quad \forall u \in H^1(\Omega). \tag{2.2}$$

Using this stronger condition on nudging strategy, we deduce that, for any $z \in V_H$,

$$\langle L_H^* z, v \rangle = (z, L_H v)_Y \leq \|z\|_Y \|L_H v\|_Y \leq C_S^{1/2} \|z\|_Y \|v\|_0 \tag{2.3}$$

and consequently, $L_H^* z \in L^2(\Omega)$. At the same time, the adjoint operator $L_H^*$ is stable as well, $\|L_H^* z\|_0^2 \leq C_S \|z\|_Y^2$ for any $z \in V_H$, since

$$\|L_H^* z\|_0^2 = (z, L_H L_H^* z)_Y \leq \frac{C_S}{2}\|z\|_Y^2 + \frac{C_S^{-1}}{2}\|L_H L_H^* z\|_Y^2 \leq \frac{C_S}{2}\|z\|_Y^2 + \frac{1}{2}\|L_H^* z\|_0^2. \tag{2.4}$$

Out of three nudging strategies presented in Section 1.4, the $L^2$ projection (1.14) and the mean value (1.18) satisfy the stability Assumption 2.2 trivially. In contrast, the boundary projection is stable, $\|P_H^\partial u\|_\partial \leq \|u\|_\partial \lesssim \|u\|_{H^{1/2}(\Omega)}$, in $H^{1/2}(\Omega)$ norm only. Nevertheless, the error analysis in Proposition 3.2 covers the nudging by the boundary projection as well, and numerical experiments of Section 4.4 demonstrate optimal rates for problems of minimal regularity.

Note that the constant $C_S$ is hidden in the rest of the paper. In this section, the symbol $\lesssim$ also hides other constants that stem from the domain $\Omega$.

### 2.1. Adjoint nudged projection

Later we will employ the Aubin–Nitsche trick by considering the solution $\psi \in H^1(\Omega)$ of the following continuous problem for a given $e \in L^2(\Omega)$,

$$a(v, \psi) + \mu(L_H v, L_H \psi)_Y = (e, v), \qquad \forall v \in H^1(\Omega). \tag{2.5}$$

**Proposition 2.2.** *For any $e \in L^2(\Omega)$, there exists a unique solution $\psi$ of (2.5) and*

$$\|\psi\|_{H^2} \lesssim (1 + \mu C_a)\|e\|_0.$$

*Proof.* Observability property (1.1) guarantees existence and uniqueness of $\psi$ by Lax–Milgram theorem as the left-hand side of (2.5) is coercive on $H^1(\Omega)$. Also, writing (1.1) for $u = \psi$ and testing (2.5) with $v = \psi \in H^1(\Omega)$ yields

$$\|\psi\|_0^2 \lesssim \|L_H \psi\|_Y^2 + C_a |\psi|_a^2 \leq \max(\mu^{-1}, C_a)(e, \psi) \leq \max(\mu^{-1}, C_a)\|e\|_0 \|\psi\|_0. \tag{2.6}$$

By writing $a(v, \psi) = (e - \mu L_H^* L_H \psi, v)$ and employing the elliptic regularity of its left-hand side, one shows that $\psi \in H^2(\Omega)$,

$$\|\psi\|_{H^2} \lesssim \|e - \mu L_H^* L_H \psi\|_0 \lesssim \|e\|_0 + \mu \|\psi\|_0 \leq (1 + \max(\mu C_a, 1))\|e\|_0 \lesssim (1 + \mu C_a)\|e\|_0,$$

where the stability (2.4) and (2.6) are used to complete the proof. $\square$

### 2.2. Regularity estimates and exponential decay

Consider the error $w = u - \nu$ of nudging and the equation it satisfies for all $t \in (0, T)$,

$$\partial_t w + (\Delta_a + \mu L_H^* L_H) w = 0, \qquad \partial_n w = 0, \qquad w(0) = -\nu_0, \tag{2.7}$$

which is obtained by subtracting the strong form of (1.1) from the nudged problem (2.1) with the observations $\nu_H$ given by (1.5). Accordingly, the weak formulation of (2.7) reads: find $w \in L^2(0, T; H^1(\Omega)) \cap H^1(0, T; H^{-1}(\Omega))$ such that, for a.e. $t \in (0, T)$,

$$\langle \partial_t w, v \rangle + a(w, v) + \mu(L_H w, L_H v)_Y = 0, \quad \forall v \in H^1(\Omega). \tag{2.8}$$

Clearly, the error equation (2.7) does not involve $f$ and $g$, and the regularity of its solution $w$ is connected to the initial data $\nu_0$. At the same time, the nudging term is sufficiently regular to guarantee that $w \in H^2(\Omega)$ for a smooth $\nu_0$, thanks to Assumption 2.1. This claim and further estimates are demonstrated in the following statement which emphasizes the dependence on $\mu$ that will be taken into account in Theorem 3.10.

**Proposition 2.3.** *The problem (2.8) is well-posed and for the error $w = u - \nu$ we have*

$$\|w\|_{L^\infty(L^2)}^2 + \mu\|L_H w\|_{L^2(Y)}^2 + \|w\|_{L^2(H^1)}^2 + (1+\mu)^{-1}\|\partial_t w\|_{L^2(H^{-1})}^2 \lesssim 1 \tag{2.9}$$

*for any $\nu_0 \in L^2(\Omega)$ (even without Assumption 2.1). Moreover, if $\nu_0 \in H^1(\Omega)$, then*

$$\|w\|_{L^\infty(H^1)}^2 + \mu\|L_H w\|_{L^\infty(Y)}^2 + \|w\|_{L^2(H^2)}^2 + \|\partial_t w\|_{L^2(L^2)}^2 \lesssim 1 + \mu\,. \tag{2.10}$$

*Proof.* The proof is standard for parabolic problems. □

Next we show that the error of nudging decays with saturated rate (1.11) thanks to the observability property (1.10) of the abstract operator $L_H$.

**Proposition 2.4** (saturated decay)**.** *Let $\mu \geq \mu_s$ in (1.4). Then, the nudged solution is guaranteed to converge exponentially in time (even without Assumption 2.1) with the rate $\gamma_s$ (1.11),*

$$\|u(t) - \nu(t)\|_0 \leq e^{-t\gamma_s}\|\nu_0\|_0\,, \qquad \forall t \in (0, T)\,. \tag{2.11}$$

*Proof.* To show the exponential decay, scale the squared norm $\|w\|_0^2$ with a yet-to-be-determined exponential factor $e^{\gamma t}$, test (2.8) with $v = w$ to obtain

$$\frac{1}{2}\frac{d}{dt}(e^{2\gamma t}\|w\|_0^2) = e^{2\gamma t}\left(-|w|_a^2 - \mu\|L_H w\|_Y^2 + \gamma\|w\|_0^2\right)\,,$$

which is non-positive if $\gamma \leq \min(\mu, 1/C_a) = \min(\mu_s, 1/C_a) = \gamma_s$ since $\|w\|_0^2 \leq \|L_H w\|_Y^2 + C_a|w|_a^2$ according to (1.10). The proof is completed by integrating over $(0, t)$, $t < T$, and recalling the initial condition in (2.7). □

*Remark* 2.5. It is clear from the proof of Proposition 2.4 that the rate of the exponential decay of $u(t)$ to $\nu(t)$ in time depends on the time-dependent value $C_a[w(t)] < C_a$ of the corresponding constant in the observability property (1.10) for $w(t) = u(t) - \nu(t)$, see Remark 1.2. Consequently, a sharper version of Proposition 2.4 holds with a time-dependent rate $\gamma(t) = (C_a[u(t) - \nu(t)])^{-1}$ and the saturated rate $\gamma_s = C_a^{-1}$ provides a global bound from below for $\gamma(t)$. This observation is relevant to the numerical examples in Section 4 when the exponential rate appear to change over time.

## 3 Error analysis

Let the space $H^1(\Omega)$ be discretized by a conforming finite-dimensional space $V_h \subset H^1(\Omega)$ based on quasi-uniform and shape-regular triangulation $\mathcal{T}_h$ of $\Omega$. The symbol $\lesssim$ hides the constants stemming from the domain $\Omega$ and from the mesh quality of $\mathcal{T}_h$. We emphasize that, if $V_H$ is a finite element space, the corresponding coarse triangulation $\mathcal{T}_H$ is not required to satisfy any properties in addition to (1.10).

In this section we assess the convergence of the following general nudging scheme on $V_h$,

$$\langle \partial_t u_h, v_h \rangle + a(u_h, v_h) + \mu(L_H u_h - \nu_H, L_H v_h)_Y = \langle f, v_h \rangle + \langle g, v_h \rangle_\partial\,, \quad \forall v_h \in V_h\,, \tag{3.1}$$

with the initial condition $u_h(0) = 0$, to the true solution $\nu$ of (1.1) with $\nu(0) = \nu_0$. The analysis utilizes the $L^2$-projection $\Pi_h : L^2(\Omega) \to V_h$ defined for a $v \in L^2(\Omega)$ by

$$(\Pi_h v, v_h) = (v, v_h)\,, \quad \forall v_h \in V_h\,. \tag{3.2}$$

**Proposition 3.1.** *The following error estimates for $L^2$-projection (3.2) hold,*

$$\|\Pi_h v - v\|_1 + h^{-1}\|\Pi_h v - v\|_0 \lesssim h^s \|v\|_{H^{1+s}}, \quad \forall v \in H^{1+s}(\Omega), \quad \forall s \geq 0. \tag{3.3}$$

*Proof.* The proof is standard an omitted therefore. $\square$

We start with the error analysis of the general nudging strategy (3.1) that follows standard arguments typically utilized for regular problems. One considers the error $u_h - \nu$, deduces the consistency equation, takes into account the approximation properties of the coarse space $V_H$ and applies Grönwall lemma to the result. However, the proof yields suboptimal rate $h^s$ for the long-term discretization error. Consequently, the resulting values $u_h(t)$ may be useless for restarting (1.8) as it can be polluted with the error $h^s$. At the same time, the proof does not require elliptic regularity of the bilinear form $a(u,v)$, i.e. the conductivity matrix does not have to be uniformly Lipschitz, and it also allows for a nudging strategy that does not satisfy stability Assumption 2.1 - both aspects are explored in the numerical experiments of Section 4.

**Proposition 3.2.** *Assume $\nu \in L^\infty(0,T; H^{1+s}(\Omega))$, $s \geq 0$. For any $\mu > 0$ in (3.1), we have*

$$\|u_h(t) - \nu(t)\|_0 \lesssim e^{-\gamma t}\|\nu_0\|_0 + h^s(1 + \mu C_\mu)^{1/2}\|\nu\|_{L^2(H^{1+s})} + h^{1+s}\|\nu(t)\|_{H^{1+s}} \tag{3.4}$$

*with $\gamma = \min(\mu, C_a^{-1})/2 \leq \gamma_s/2$ and $C_\mu \lesssim 1$. In addition, $C_\mu \simeq h^2$ if Assumption 2.1 holds.*

*Proof.* We utilize the $L^2$-projection $\Pi_h$ by splitting the error of nudging as follows,

$$u_h - \nu = \zeta_h + e_\nu^\Pi, \quad \zeta_h := u_h - \Pi_h \nu \in V_h, \quad e_\nu^\Pi := \Pi_h \nu - \nu \in H^1(\Omega).$$

One express the control of discrete error $\zeta_h$ by the consistency equation for (3.1),

$$(\partial_t \zeta_h, v_h) + a(\zeta_h, v_h) + \mu(L_H \zeta_h, L_H v_h)_Y = (f, v_h) + (g, v_h)_\partial + \mu(\nu_H, L_H v_h)_Y \tag{3.5}$$
$$- (\partial_t \Pi_h \nu, v_h) - a(\Pi_h \nu, v_h) - \mu(L_H \Pi_h \nu, L_H v_h)_Y = -a(e_\nu^\Pi, v_h) - \mu(L_H e_\nu^\Pi, L_H v_h)_Y,$$

where for the last equality we utilized (1.5) and that $(\partial_t \Pi_h \nu, v_h) = (\partial_t \nu, v_h)$ thanks to (3.2).

To show the exponential decay under the observability condition (1.10), we test (3.5) with $v_h = \zeta_h$, scale $\|\zeta_h\|_0^2$ with a yet-to-be-determined exponential factor $e^{\gamma t}$ and obtain that

$$\frac{1}{2}\frac{d}{dt}(e^{2\gamma t}\|\zeta_h\|_0^2) = e^{2\gamma t}\left(\gamma\|\zeta_h\|_0^2 - |\zeta_h|_a^2 - \mu\|L_H\zeta_h\|_Y^2 - a(e_\nu^\Pi, \zeta_h) - \mu(L_H e_\nu^\Pi, L_H \zeta_h)_Y\right)$$
$$\leq e^{2\gamma t}\left(\gamma\|\zeta_h\|_0^2 - \frac{1}{2}|\zeta_h|_a^2 - \frac{1}{2}\mu\|L_H\zeta_h\|_Y^2 + \frac{1}{2}|e_\nu^\Pi|_a^2 + \frac{1}{2}\mu\|L_H e_\nu^\Pi\|_Y^2\right).$$

The observability property, $\|\zeta_h\|_0^2 \leq \|L_H\zeta_h\|_Y^2 + C_a|\zeta_h|_a^2$, the choice $\gamma = \frac{1}{2}\min(\mu, 1/C_a)$ and the integration in time imply the estimate

$$\|\zeta_h(t)\|_0^2 \leq e^{-2\gamma t}\|\zeta_h(0)\|_0^2 + \int_0^t \left(\frac{1}{2}|e_\nu^\Pi|_a^2 + \frac{1}{2}\mu\|L_H e_\nu^\Pi\|_Y^2\right).$$

Noticing $\|\zeta_h(0)\|_0 = \|\Pi_h \nu_0\|_0 \leq \|\nu_0\|_0$, estimating by (1.6) (or by (2.2) if Assumption 2.1 holds) and applying triangle inequality to $\zeta_h + e_\nu^\Pi$ followed by Proposition 3.1 complete the proof. $\square$

In the spirit of [DLR25], the next corollary highlights that the dependence on $\mu$ of the error estimates above disappears if the projection $\Pi_h$ and operator $L_H$ are related as specified in the statement. The strategy (1.18) of nudging by the mean value is an example.

**Corollary 3.3.** *Let $\infty > \mu \geq \mu_s$ in Proposition 3.2. If $L_H \Pi_h \nu(t) = L_H \nu(t)$ in $V_H$, then*

$$\|u_h(t) - \nu(t)\|_0 \lesssim e^{-t/(2C_a)} + h^s. \tag{3.6}$$

*Proof.* For a large $\mu$ we still have $\gamma = \gamma_s/2$ in (3.4). Also, since the coarse data $\nu_H$ is given by $L_H\nu = L_H\Pi_h\nu$, the term $L_H e_\nu^\Pi$ vanishes at the end of the proof above. $\square$

## 3.1. Nudged elliptic projection

Before discussing problems of limited regularity, we introduce the nudged elliptic projection $Q_h$ defined in (3.7) that improves Proposition 3.2 by yielding optimal rates for the long-term discretization error in the $L^\infty(0,T,L^2(\Omega))$ sense under some additional regularity in time. As the result, the values $u_h(t)$ can be used in restarting of (1.8) in this case.

Consider the following discrete minimization problem: Given $u \in H^1(\Omega)$, find $Q_h u \in V_h$ such that
$$N_\mu[u_h] = \frac{1}{2}|u_h - u|_a^2 + \frac{1}{2}\mu \|L_H(u_h - u)\|_Y^2$$
is minimized over $V_h$. The weak formulation for the minimizer $Q_h u$ is obtained by taking a discrete variation of $N_\mu$ and setting it to zero as in
$$a(Q_h u, v_h) + \mu(L_H Q_h u, L_H v_h)_Y = a(u, v_h) + \mu(L_H u, L_H v_h)_Y, \quad \forall v_h \in V_h. \tag{3.7}$$

The minimization problem (3.7) is well-posed thanks to observability property (1.1). Moreover, the minimization error $e_u^Q := Q_h u - u$ possesses Galerkin orthogonality and
$$a(e_u^Q, v_h) + \mu(L_H e_u^Q, L_H v_h)_Y = 0, \quad N_\mu[Q_h u] \leq N_\mu[v_h], \quad \forall v_h \in V_h. \tag{3.8}$$

Note that a similar to (3.7) construction appears in [DLR25]. We now derive error estimates for the nudged projection $Q_h u$.

**Lemma 3.4** (nudged elliptic projection). *For $u \in H^{1+s}(\Omega)$, $s \geq 0$, we have*
$$|Q_h u - u|_a + \mu^{1/2}\|L_H(Q_h u - u)\|_Y + (h r_h q_H)^{-1}\|Q_h u - u\|_0 \lesssim h^s r_h |u|_{H^{1+s}} \tag{3.9}$$
*with the following constants that depend on $h$, $\mu$ and $C_a$ from (1.10),*
$$r_h = 1 + h\mu^{1/2}, \quad q_H = 1 + \mu C_a. \tag{3.10}$$

*Proof.* Best approximation property (3.8) for $v_h = \Pi_h u \in H^1(\Omega)$ and estimates from Proposition 3.1 yield for $e_u^Q := Q_h u - u$ that
$$|e_u^Q|_a^2 + \mu \|L_H e_u^Q\|_Y^2 \leq |\Pi_h u - u|_a^2 + \mu \|L_H(\Pi_h u - u)\|_Y^2 \lesssim h^{2s} r_h^2 |u|_{H^{1+s}}^2, \tag{3.11}$$
where $(1 + \mu h^2) \leq (1 + h\mu^{1/2})^2 =: r_h^2$. For an estimate of $\|e_u^Q\|_0$, consider the solution $\psi$ to the adjoint problem (2.5) for $e := e_u^Q$, test (2.5) with $v = e_u^Q$ and subtract (3.8) tested with $v_h = I_h \psi$ from the result to derive
$$\|e_u^Q\|_0^2 = a(e_u^Q, \psi) + \mu(L_H e_u^Q, L_H \psi)_Y = a(e_u^Q, \psi - \Pi_h \psi) + \mu(L_H e_u^Q, L_H(\psi - \Pi_h \psi))_Y$$
$$\lesssim |e_u^Q|_a |\psi - \Pi_h \psi|_a + \mu \|L_H e_u^Q\|_Y \|\psi - \Pi_h \psi\|_0 \lesssim h(|e_u^Q|_a + h\mu \|L_H e_u^Q\|_Y)\|\psi\|_{H^2},$$
where (2.2) and the estimates (3.1) are used in the second line. Combining the estimate from Proposition 2.2 with (3.11) gives
$$\|e_u^Q\|_0^2 \lesssim h(1 + h\mu^{1/2}) h^s r_h |u|_{H^{1+s}} \|\psi\|_{H^2} \lesssim h^{s+1} r_h^2 (1 + \mu C_a) |u|_{H^{1+s}} \|e_u^Q\|_0 \tag{3.12}$$
thus completing the proof with $q_H := 1 + \mu C_a$. □

**Proposition 3.5.** *Assume $\nu \in H^1(0, T; H^{1+s}(\Omega))$, $s \geq 0$. For any $\mu > 0$ in (3.1), we have*
$$\|u_h(t) - \nu(t)\|_0 \lesssim e^{-\gamma t}\|\nu_0\|_0 + h^{1+s} r_h^2 q_H(\|\partial_t \nu\|_{L^2(H^{1+s})} + |\nu(t)|_{H^{1+s}})$$
*with $0 < \gamma + 1 = \min(\mu, C_a^{-1}) \leq \gamma_s$ and $r_h$, $q_H$ given in (3.10)*

*Proof.* One repeats the arguments for Proposition 3.2 with minor adjustments which are now outlined briefly. Introducing the error $e_\nu = u_h - \nu = (u_h - Q_h\nu) + (Q_h\nu - \nu) =: \theta_h + e_\nu^Q$, we deduce the following control of $\theta_h$,

$$(\partial_t \theta_h, v_h) + a(\theta_h, v_h) + \mu(L_H \theta_h, L_H v_h)_Y = -\left(\partial_t e_\nu^Q, v_h\right), \tag{3.13}$$

using Galerkin orthogonality (3.8) of the nudged elliptic projection $Q_h u$. Estimating the term on the right of (3.13) by Young's inequality, applying observability property (1.10) for $\theta_h$, noticing $\theta_h(0) = Q_h \nu_0$ yield upon integration in time the final estimate,

$$\|e_\nu(t)\|_0^2 \le e^{-2\gamma t}\|Q_h\nu_0\|_0^2 + \|e_\nu^Q(t)\|_0^2 + \frac{1}{4}\int_0^t \|\partial_t e_\nu^Q\|_0^2.$$

Application of Lemma 3.4 completes the proof. $\square$

*Remark* 3.6. The estimate of Proposition 3.5 depends on $\mu$ only through the constants $r_h$ and $q_H$. This is in contrast to the estimate of Proposition 3.2 where the dependence on $\mu$ can be traced back to (3.5). Therefore, potential improvement of Lemma 3.4 in this regard leads to $\mu$-independent error estimates for the nudging scheme (3.1).

The next corollary of Proposition 3.5 corresponds to the saturated choice (1.12) of the nudging parameter for which the resulting constant $r_h^2 = 1 + h^2 \mu_s = 1 + h^2/C_a$ solely causes the dependence on on $C_a$ of the long-term discretization error.

**Corollary 3.7.** *Let $\mu = \mu_s$ in Proposition 3.5, then*

$$\|u_h(t) - \nu(t)\|_0 \lesssim e^{-t/C_a} + h^{1+s}(1 + h^2/C_a).$$

## 3.2. Limited regularity

Before this section, we have been studying the error $\nu - u_h$. Now we analyze directly the following discretization error of the nudged problem (1.4),

$$e_u = u_h - u, \tag{3.14}$$

assuming minimally that the nudged solution $u \in L^2(0, T; H^{1+s}(\Omega))$, $s \ge 0$, supported by the regularity analysis of the nudged problem in Proposition 2.3. Lemma 3.15 shows the stability of the scheme (3.1) and, following [CH02], we employ the parabolic version of the duality argument in Lemma 3.9 to culminate in main Theorem 3.10.

**Lemma 3.8.** *Assume $u \in L^2(0, T; H^{1+s}(\Omega))$, $s > 0$. For any $\mu > 0$ we have*

$$\int_0^t |e_u|_a^2 + \mu \int_0^t \|L_H e_u\|_Y^2 \lesssim h^{2s}(1 + \mu h^2) \int_0^t \|u\|_{H^{1+s}}^2. \tag{3.15}$$

*Proof.* We utilize the $L^2$-projection $\Pi_h$ by splitting the error (3.14) as follows,

$$u_h - u = \theta_h + e_u^\Pi, \qquad \theta_h := u_h - \Pi_h u \in V_h, \quad e_u^\Pi := \Pi_h u - u \in H^1(\Omega).$$

We deduce the control of the discrete error $\theta_h$ as in (3.5),

$$(\partial_t \theta_h, v_h) + a(\theta_h, v_h) + \mu(L_H \theta_h, L_H v_h)_Y = -a(e_u^\Pi, v_h) - \mu(L_H e_u^\Pi, L_H v_h)_Y. \tag{3.16}$$

Testing (3.16) with $v_h = \theta_h$, integrating in time and using Cauchy–Schwarz yield

$$\|\theta_h(t)\|_0^2 + \int_0^t |\theta_h|_a^2 + \mu \int_0^t \|L_H \theta_h\|_Y^2 \lesssim \int_0^t |e_u^\Pi|_a^2 + \mu \int_0^t \|L_H e_u^\Pi\|_Y^2 \qquad (3.17)$$

taking into account $\theta_h(0) = u_h(0) = \Pi_h u_0 = 0$. Therefore, we have

$$\int_0^t |\theta_h + e_u^\Pi|_a^2 + \mu \int_0^t \|L_H(\theta_h + e_u^\Pi)\|_Y^2 \lesssim h^{2s}(1 + \mu h^2) \int_0^t \|u\|_{H^{1+s}}^2$$

thanks to the error estimates for $\Pi_h$-projection from Lemma 3.1. $\square$

The forthcoming Lemma 3.9 utilizes the parabolic duality argument for the error $e_u = u_h - u$ and restores the optimality, see Remark 3.11.

**Lemma 3.9.** *For any $\mu > 0$ we have the following estimate with $r_h$ and $q_H$ given in (3.10),*

$$\int_0^t \|e_u\|_0^2 \lesssim h^2 r_h^2 q_H^2 \left( \int_0^t |e_u|_a^2 + \int_0^t \mu \|L_H e_u\|_{L^2(Y)}^2 \right). \qquad (3.18)$$

*Proof.* Recall $e_u = u_h - u \in L^2(0, T; H^1(\Omega))$ and note that $e_u(0) = u(0) = u_h(0) = 0$. From Proposition 2.2, there exists $\psi(t) \in H^2(\Omega)$ that solves the adjoint nudged problem (2.5) with $e := e_u(t)$ a.e. $t \in (0, T]$ such that $\|\psi(t)\|_{H^2} \lesssim q_H \|e_u(t)\|_0$ with $q_H$ defined in (3.10).

Now consider the nudged projection $\psi_h(t) := Q_h \psi(t)$ a.e. $t \in (0, T]$ and deduce from (2.5) and (3.7) that it satisfies the discrete problem $a(v_h, \psi_h(t)) + \mu(L_H v_h, L_H \psi_h(t))_Y = (e_u(t), v_h)$, for all $v_h \in V_h$, and therefore, a.e. $t \in (0, T]$,

$$a(\partial_t \psi_h, \psi_h) + \mu(L_H(\partial_t \psi_h), L_H \psi_h)_Y = \langle \partial_t e_u, \psi_h \rangle. \qquad (3.19)$$

In addition, from Lemma 3.4, we have the following estimate of the error $e_\psi^Q := \psi - Q_h \psi$,

$$|e_\psi^Q|_a^2 + \mu \|L_H e_\psi^Q\|_Y^2 \lesssim h^2 r_h^2 |\psi|_{H^2}^2 \lesssim h^2 r_h^2 q_H^2 |e_u|_0^2. \qquad (3.20)$$

The error $e_u = u_h - u$ satisfies the discretized error equation (2.8) for all $t \in (0, T)$,

$$(\partial_t e_u, v_h) + a(e_u, v_h) + \mu(L_H e_u, L_H v_h)_Y = 0, \quad \forall v_h \in V_h.$$

Taking $v_h = \psi_h$ in the above and testing (2.5) with $v = e_u$ yield

$$\|e_u\|_0^2 = a(\psi, e_u) + \mu(L_H \psi, L_H e_u)_Y = a(e_\psi^Q, e_u) + \mu(L_H e_\psi^Q, L_H e_u)_Y - \langle \partial_t e_u, \psi_h \rangle. \qquad (3.21)$$

The first part of (3.21) is estimated by utilizing (3.20) as follows,

$$a(e_\psi^Q, e_u) + \mu(L_H e_\psi^Q, L_H e_u)_Y \leq (|e_\psi^Q|_a^2 + \mu \|L_H e_\psi^Q\|_Y^2)^{1/2} (|e_u|_a^2 + \mu \|L_H e_u\|_Y^2)^{1/2}$$

$$\leq h r_h q_H \|e_u\|_0 (|e_u|_a^2 + \mu \|L_H e_u\|_Y^2)^{1/2} \leq \frac{1}{2} \|e_u\|_0^2 + \frac{1}{2} h^2 r_h^2 q_H^2 (|e_u|_a^2 + \mu \|L_H e_u\|_Y^2).$$

For the remaining part in (3.21), recall (3.19) and integrate it in time from 0 to $t$,

$$2 \int_0^t \langle -\partial_t e_u, \psi_h \rangle = \int_0^t -\frac{d}{dt}(a(\psi_h, \psi_h) + \mu \|L_H \psi_h\|_Y^2) = |\psi_h(0)|_a^2 + \mu \|L_H \psi_h(0)\|_Y^2$$

$$- |\psi_h(t)|_a^2 - \mu \|L_H \psi_h(t)\|_Y^2 \leq |\psi_h(0)|_a^2 + \mu \|L_H \psi_h(0)\|_Y^2 = (e_u(0), \psi_h(0)) = 0.$$

Consequently, one integrates (3.21) over $(0, t)$ to find

$$\int_0^t \|e_u\|_0^2 \lesssim \frac{1}{2} \int_0^t \|e_u\|_0^2 + \frac{1}{2} h^2 r_h^2 q_H^2 \left( \int_0^t |e_u|_a^2 + \int_0^t \mu \|L_H e_u\|_{L^2(Y)}^2 \right),$$

and the first term on the right-hand side is absorbed to yield the claim. $\square$

We are in position to prove the main result below in which the nudging parameter is saturated, $\mu = \mu_s$, for the clarity of presentation only. Note that the time interval $(t, T)$ is involved in Theorem 3.10 which guarantees that the nudged solution $u_h$ is close to $\nu$ on $(t, T)$. The corresponding norm and the long-term discretization error are optimal taking into account the assumed regularity; see Remark 3.11. Note that it implies that we cannot use $u_h(t)$ as a restarting initial condition in (1.8) and the nudging scheme (3.1) must be run for all $t \in (0, T)$.

**Theorem 3.10.** *Assume $\nu \in L^2(0, T; H^{1+s}(\Omega))$, $0 \leq s \leq 1$, solves (1.1) and $\nu_0 \in H^1(\Omega)$ in (1.3). For the saturated choice $\mu = \mu_s$ (1.12) in (3.1) and with $C_a$ defined in (1.10), we have*

$$\|u_h - \nu\|_{L^2(t,T;L^2(\Omega))} \lesssim C_a^{1/2} e^{-t/C_a} + h^{s+1}(1 + h^2/C_a)(1 + 1/C_a)^{1/2}. \tag{3.22}$$

*Proof.* Combining Lemma 3.8 and Lemma 3.9 with $\mu = \mu_s = C_a^{-1}$, one derives

$$\|u_h - u\|_{L^2(t,T;L^2(\Omega))} \lesssim \|u_h - u\|_{L^2(0,T;L^2(\Omega))} \lesssim h^{1+s}(1 + h^2 \mu_s)(1 + \mu_s)^{1/2}. \tag{3.23}$$

Since Proposition 2.3 implies that $u = \nu + (u - \nu) \in L^2(0, T; H^{1+s}(\Omega))$ and

$$\int_t^T \|u\|_{H^{1+s}}^2 \lesssim \int_t^T \|\nu\|_{H^{1+s}}^2 + \int_t^T \|u - \nu\|_{H^2}^2 \lesssim 1 + \mu = 1 + \mu_s.$$

On another hand, we utilize the estimate from Proposition 2.4 by integrating it over $(t, T)$,

$$\|u - \nu\|_{L^2(t,T;L^2(\Omega))}^2 \lesssim \int_t^T e^{-2t\gamma_s} \|\nu_0\|_0^2 \lesssim C_a e^{-2t/C_a}(1 - e^{-2(T-t)/C_a}). \tag{3.24}$$

Triangle inequality in $u_h - \nu = (u_h - u) + (u - \nu)$ completes the proof. $\square$

*Remark* 3.11. One expects an error estimate $\|\nu_h - \nu\|_{L^2(0,T;L^2(\Omega))} \lesssim h^{s+1}$ for the discrete problem (1.8) with a known initial condition in the case of limited regularity, $\nu \in L^2(0, T; H^{1+s}(\Omega))$. Consequently, the results in Theorem 3.10 are optimal for a fixed observation space $V_H$ and, therefore, a fixed $C_a$. Alternatively, if one does not seek to maximize the exponential rate of convergence in time by saturating the nudging parameter, we can set $\mu$ in (3.1) to any constant $\mu_0 < \mu_s$ and Lemma 3.8 followed by Lemma 3.9 imply similarly to Theorem 3.10 that

$$\|u_h - \nu\|_{L^2(t,T;L^2(\Omega))} \lesssim e^{-t\mu_0} + h^{s+1}(1 + \mu_0) \tag{3.25}$$

## 4 Numerical experiments

In this section we conduct a detailed numerical investigation of the nudging strategies from Section 1.4 and support the semi-discrete error analysis presented in Section 3. To this end, each strategy is applied for three test problems introduced in Section 4.1 as they differ in several key aspects including the regularity of solutions and the continuity of the conductivity coefficients. The common discretization details can be found in Section 4.2 including the rationale of the numerical experiments – saturation and convergence tests. The numerical results are grouped by the nudging strategy involved: Section 4.3 for the nudging by the approximating projection operator, Section 4.4 for the nudging by boundary projection, and Section 4.5 for the nudging by the mean value.

### 4.1. Test problems

Here we specify three exact parabolic problems (1.1) in the domain $\Omega = [-1,1]^2$ with $T = 3$ that are later involved in the numerical experiments. The boundary data is computed explicitly as $g = \nabla \nu \cdot \mathbf{n}$ using the outer normal $\mathbf{n}$ to $\partial \Omega$ for each test case. These manufactured solutions are either time-independent with $\omega = 0$ for saturation tests or time-dependent with $\omega = \pi$ for the convergence tests.

- **A smooth problem.** Consider the solution $\nu \in C^\infty(0,T;C^\infty(\Omega))$ of (1.1) with
$$a(\nu, v) = (\nabla \nu, \nabla v), \qquad f = -\omega \sin(\omega t)|\mathbf{x} - \mathbf{x}_0|^2 - 4\cos(\pi t), \tag{4.1}$$
where $\mathbf{x}_0 = (0.5, 0.5) \in \Omega$, given by
$$\nu(\mathbf{x}, t) = \cos(\omega t)|\mathbf{x} - \mathbf{x}_0|^2. \tag{4.2}$$

- **Dirac delta source.** Consider the solution $\nu \in C^\infty(0, T; H^1(\Omega))$ of (1.1) with
$$a(u, v) = (\nabla u, \nabla v), \qquad f = \cos(\pi t)\delta(\mathbf{x}_0) - \frac{\omega \sin(\omega t)}{2\pi}\ln(|\mathbf{x} - \mathbf{x}_0|), \tag{4.3}$$
where $\delta(\mathbf{x}_0)$ is the Dirac delta located at the singularity $\mathbf{x}_0 = (1/3, 1/3)$, given by
$$\nu(\mathbf{x}, t) = \frac{\cos(\omega t)}{2\pi}\ln(|\mathbf{x} - \mathbf{x}_0|). \tag{4.4}$$

Although it is not the focus of the paper, we would like to mention that the accuracy of the numerical quadratures turns out to be critical for performance of the nudged scheme.

- **Kellogg problem with discontinuous conductivity [Bru74; BDN13].** Consider the solution $\nu \in C^\infty(0, T; H^{5/4-\epsilon}(\Omega))$, for any $\epsilon > 0$, of (1.1) with
$$a(u, v) = (A(x,y)\nabla u, \nabla v), \qquad f = -\omega \sin(\omega t) N(\mathbf{x}). \tag{4.5}$$

Here $N(\mathbf{x})$ is given in (4.7) and the conductivity matrix $A(x,y)$ is proportional to the identity matrix $I$ with the following discontinuous coefficients,
$$A(x,y) = \begin{cases} bI, & \text{if } xy \geq 0 \\ I, & \text{otherwise}. \end{cases} \tag{4.6}$$

The exact solution is given in the polar coordinates $\mathbf{x} = (\rho, \theta)$ by
$$\nu(\mathbf{x}, t) = \cos(\omega t) N(\mathbf{x}), \qquad N(\mathbf{x}) = N(\rho, \theta) = \rho^\alpha \mu(\theta), \tag{4.7}$$
where $\mu(\theta)$ is a specific function which has kinks at $\theta = 0, \frac{\pi}{2}, \pi, \frac{3\pi}{2}$ to compensate the discontinuity lines $x = 0$ and $y = 0$ of $A$. The details can be found in [BDN13, Section 6.2] and we only mention the values of the following key parameters,
$$\alpha = 0.25, \qquad b = 25.27414236908818, \tag{4.8}$$
which guarantee that $N(\mathbf{x}) \in H^{1+s}(\Omega)$ in (4.7) for any $s < \alpha$. Note that the meshes in all numerical test are aligned with the lines of discontinuity of $A$ for accuracy reasons. We emphasize that $\mu(\theta)$ taken from [BDN13] has the zero mean value and, consequently, the nudging by the mean value becomes useless as the nudged scheme is initialized with $u_0 = 0$, $\bar{u}_0 = 0$. At the same time, since the conductivity matrix (4.6) is discontinuous, the elliptic regularity assumed throughout the paper does not hold for the Kellogg problem. In this case, Theorem 3.10 and Proposition 3.5 are not applicable but Proposition 3.2 still guarantees convergence with rates.

## 4.2. Discrete setting of tests

For comparison reasons, the same discrete setting is used in all numerical experiments and is described hereafter. The observational data of all nudging strategies considered are based on the uniform quadrilateral triangulation $\mathcal{T}_H$ of $\Omega = [-1,1]^2$ with the following *coarse* mesh size,

$$H = 2^{-L}, \qquad L = 2, \tag{4.9}$$

except for the mean value nudging that requires no triangulation. For the discrete problems (3.1) we choose a sequence of uniform quadrilateral triangulations $\mathcal{T}_h$ of $\Omega$ given in the following space and time refinement strategies,

$$h = 2^{-\ell}, \quad \tau = T/N, \quad N = 2^{2\ell-3}, \qquad \ell = 2, ..., 7. \tag{4.10}$$

Here $\tau$ is the time step of the implicit Euler method for $n = 0, .., N - 1$ that is initialized with $u_h^0 = 0$ and that discretizes the time derivative as in (3.1) by

$$\partial_\tau u_h := (u_h^{n+1} - u_h^n)/\tau, \qquad V_h := Q_1(\mathcal{T}_h), \tag{4.11}$$

where $Q_1(\mathcal{T}_h)$ is the conforming Lagrange finite element space of degree 1 associated with the triangulation $\mathcal{T}_h$.

For reference, we present the error $\|\nu_h(t_k) - \nu(t_k)\|_0$, $k = 0, ..., N$, of the fully discretized scheme (1.8) using $\ell = 6$ and $\ell = 7$ in (4.10) for the tests of Section 4.2.1 and Section 4.2.2, correspondingly. To compare with the nudged solutions with $u_h(0) = 0$, we plot the evolution of (1.8) with the *correct* initial condition, $\nu_h(0) = \Pi_h \nu_0$, and with the *zero* initial condition, $\nu_h(0) = 0$.

### 4.2.1. Saturation tests

For each nudging strategy, we study the saturation effect (1.11) in every test case of Section 4.1 with $\omega = 0$ by fixing the refinement level $\ell = 6$ in (4.10), setting the increasing values of the nudging parameter $\mu$ in (3.1) and evaluating the evolution of the error of nudging in time,

$$\|u_h(t_k) - \nu(t_k)\|_0, \qquad 0 \leq t_k = k\tau \leq T, \qquad k = 0, ..., N. \tag{4.12}$$

Moreover, we provide with the experimental estimation of the saturated rate $\gamma$ of the exponential decay guaranteed by Proposition 2.4; see Figure 1 for the nudging by FE projection, Figure 3 for the nudging by boundary projection, and Figure 5 for the nudging by the mean value.

The results show that the solutions to the nudged scheme (3.1) with different nudging strategies converge exponentially in time to the exact solution of (1.1). Moreover, we see that the saturation (1.11) manifests in all test problems highlighting that the range of practically useful values of $\mu$ is limited by (1.12). However, in the case of Kellogg problem, this manifestation is not clear, and we refer to Remark 1.2 for an explanation.

### 4.2.2. Convergence tests

To access the numerical accuracy of each nudging strategy in (3.1), we run the convergence experiments by increasing the refinement level $\ell$ in (4.10) while keeping a fixed value of $\mu = 64$ and plotting the evolution of the nudging error (4.12) for every test case of Section 4.1 with $\omega = \pi$; see Figure 2 for the nudging by FE projection, Figure 4 for the nudging by boundary projection, and Figure 6 for the nudging by the mean value.

In addition, we report rates of convergence with respect to the refinement level $\ell$ of the error (4.12) (along with the semi-$H^1$ norm) at the final time $t_N = T$, and of the discrete analog of the error in Theorem 3.10, i.e. the rates of the following quantities,

$$\left(\sum_{k=M}^{N} \Delta t \|u_h(t_k) - \nu(t_k)\|_0^2\right)^{1/2}, \qquad \|u_h(T) - \nu(T)\|_0, \qquad \|\nabla(u_h(T) - \nu(T))\|_0, \qquad (4.13)$$

where a sufficiently large integer $M$, $0 < M < N$, is chosen during the post-processing to guarantee that by the time $T_\mu = M\tau$ the exponential decay term in the estimate (3.22) of $\|u_h(t) - \nu(t)\|_{L_2(T_\mu, T, L^2)}$ is dominated by the discretization error.

The observed rates of convergence in different test problems with respect to the refinement (3.1) demonstrate that the nudging scheme (3.1) is as good as the non-nudged scheme (1.8) with the known initial condition after time $T_\mu = M\tau$, i.e.

$$\|u_h(t) - \nu(t)\|_{L^2(T_\mu, T, L^2)} \simeq \|\nu_h(t) - \nu(t)\|_{L^2(T_\mu, T, L^2)} \simeq h^{1+s}, \qquad (4.14)$$

where $s$ indicates the spatial regularity of the exact solution $\nu$, see Theorem 3.10.

### 4.3. Nudging by FE projection

Consider the nudging strategy (1.14) with the observations given in $V_H := Q_1(\mathcal{T}_H)$. The full discretization (4.11) of the nudging scheme (1.4) reads: Given $f^{n+1} := f(t_{n+1})$, $g^{n+1} := g(t_{n+1})$, $\nu^{n+1} := \nu(t_{n+1})$, $n = 0, ..., N-1$, find $u_h := u_h^{n+1} \in V_h$ and $z_H := z_H^{n+1} \in Q_1(\mathcal{T}_H)$ satisfying for $n = 0, ..., N-1$ the following coupled linear systems,

$$(\partial_\tau u_h, v_h) + a(u_h, v_h) + (z_H, v_h) = \langle f^{n+1}, v_h \rangle + \langle g^{n+1}, v_h \rangle_\partial, \quad \forall v_h \in V_h, \qquad (4.15)$$

$$(u_h, r_H) - \mu^{-1}(z_H, r_H) = (\nu^{n+1}, r_H), \quad \forall r_H \in Q_1(\mathcal{T}_H), \qquad (4.16)$$

where the second equation is the weak definition of the *nudger* $z_H \in Q_1(\mathcal{T}_H)$,

$$z_H = \mu P_H(u_h - \nu), \qquad \mu(L_H u_h - L_H \nu, L_H v_h)_Y = (z_H, P_H v_h) = (z_H, v_h), \qquad (4.17)$$

since $L_H := P_H$ in (1.14) is the $L^2$-projection onto $Q_1(\mathcal{T}_H)$ which also implies $(L_H u_h, r_H) = (u_h, r_H)$ and $(L_H \nu, r_H) = (\nu, r_H)$.

The nudging scheme (4.15) is subjected first to the saturation test from Section 4.2.1 for the values $\mu = \{4, 64, 1024, 16384\}$ of the nudging parameter. Figure 1 demonstrates the resulting evolutions of the error of nudging over time for different test problems, and we now comment on the numerical observations. In the case of the smooth problem (top, left), the nudging scheme (4.15) with $\mu = 4$ shows exponential decay unlike the non-nudged scheme (1.8) (dashed) with zero initial condition. Note that the evolution of (1.8) (solid) with the correct initial condition is closer to machine epsilon in magnitude for the smooth problem (4.2) and is not visible for this reason. Increasing the nudging parameter to $\mu = 64$ accelerates the nudging as its slope indicates the boost of the exponential decay. Greater values, $\mu = 1024$ and $\mu = 16384$, improve the initial stage of evolution when $t < 0.05$, but further slopes of exponential rate are similar to that of $\mu = 64$. Consequently, the saturation effect manifests and the observed limiting slope $\gamma$ is approximately depicted via a triangle and is quantified as well (bottom, right). Analogous situation can be seen for the minimally regular problem with Dirac delta (top, right). The saturation manifests around $\mu = 64$, and $\mu = 1024$ and $\mu = 16384$ coincide completely after the initial stage, i.e. for $t > 0.01$. The nudging scheme drops the error to the discretization level by $t = 0.06$ as can be seen from the evolution of the non-nudged scheme (1.8) with the correct initial condition (solid). The observed saturated rate $\gamma$ is depicted by a triangle whose slope is reported as well (bottom, right). In the case of Kellogg problem it is not clear whether the saturation manifests or not as further increase of $\mu$ keeps improving the nudging, although both $\mu = 1024$ and $\mu = 16384$ seem to nudge the solution to the discretization level by the same time $t = 0.01$. The reported observed rate is much higher than the rates in the cases of smooth and Dirac problems but we would like to remind that the saturated parameter (1.12) guarantees the rate (1.11) but does not prohibit an acceleration, see Remark 1.2. In fact, one can see similar accelerations during the initial stages for the smooth and Dirac problems.

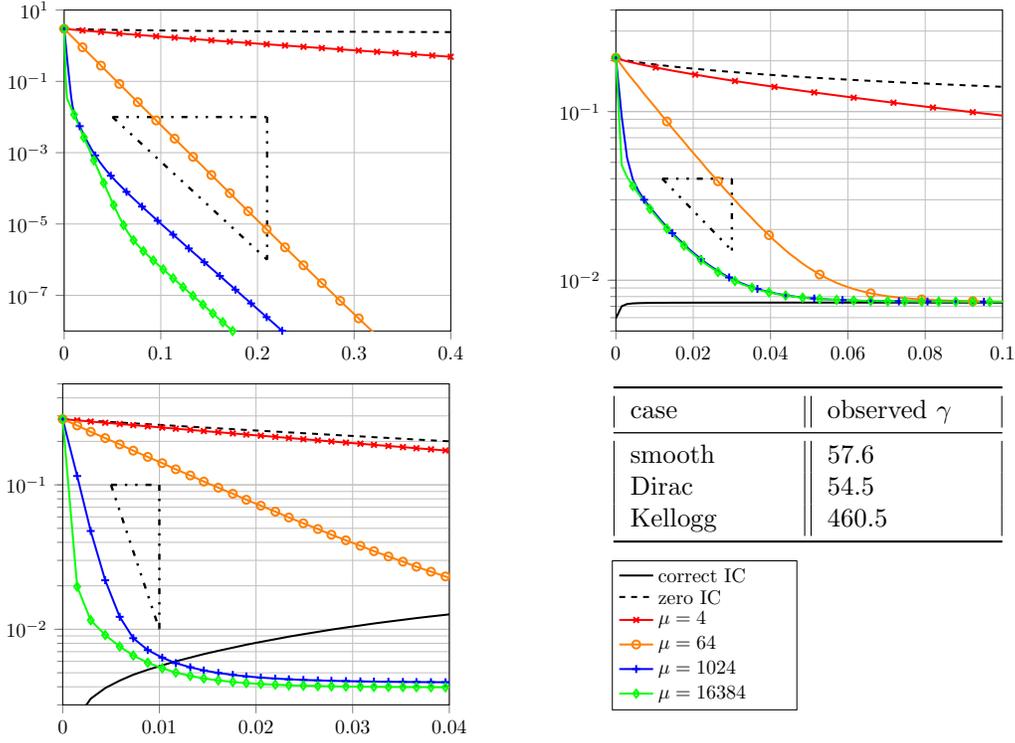

Figure 1: Saturation test described in Section 4.2.1 for the nudging by FE projection (1.14). Errors (4.12) over time for several values of $\mu$ and for the smooth (4.2) (top, left), Dirac (4.4) (top, right) and Kellogg (4.7) (bottom, left) problems. Observed rate $\gamma$ in (1.11) (bottom, right, with the common legend) correspondent to reference triangles. See Section 4.2 for discretization details and Section 4.3 for a discussion of results.

Next we run the convergence test from Section 4.2.2 for values $\ell = \{4, 5, 6, 7\}$ in (4.10). On Figure 2 we present the evolution of nudging for test problems from Section 4.1 and the rates of convergence for errors (4.13). We now briefly comment on the observations. According to Figure 1, the chosen $\mu = 64$ is sufficiently close to the saturated value (1.12) for the smooth and Dirac test problems as the observed exponential rate $\gamma < 60$ in their cases. Moreover, all nudging strategies reduce the error to the discretization level represented by the solid lines by the time $t = 0.2$. For the smooth and the Dirac problem, the rates of convergence of the nudging error at the final time $t = 3$ are one and two, correspondingly, as guaranteed by Proposition 3.5 and regularities of solutions stated in Section 4.1. Since $\mu = 64$ is sufficintly large according to saturation tests, the exponential term in Theorem 3.10 falls below the discretization level by $t = 0.2$. For this reason we report the error of nudging in $L^2(0.4, 3; L^2(\Omega))$ norm. The smooth and Dirac problem show the first and second order rates in this norm. The Kellogg problem (4.7), as mentioned before, does not possess elliptic regularity and Proposition 3.5 and Theorem 3.10 do not apply. Nevertheless, Proposition 3.2 guarantees that the rate is at least 0.25 as the Kellogg solution belongs to $H^{5/4-\varepsilon}(\Omega)$, for all $\varepsilon > 0$. In this test case, the rates of $L^2$ errors of nudging seem to double, i.e. the rate is closer to 0.5. The observed rates of convergence fully support the error analysis of Section 3.

## 4.4. Nudging by the boundary projection

If the nudging strategy corresponds to (1.16) then $V_H := Q_0(\partial \mathcal{T}_H)$ is the observational space, and the full discretization (4.11) of the nudging scheme (1.4) reads: Given $f^{n+1} := f(t_{n+1})$, $g^{n+1} := g(t_{n+1})$, $\nu^{n+1} := \nu(t_{n+1})$, $n = 0, ..., N-1$, find $u_h := u_h^{n+1} \in V_h$ and $z_H := z_H^{n+1} \in Q_0(\partial \mathcal{T}_H)$ satisfying the following coupled linear system,

$$(\partial_\tau u_h, v_h) + a(u_h, v_h) + (z_H, v_h)_\partial = \langle f^{n+1}, v_h \rangle + \langle g^{n+1}, v_h \rangle_\partial, \quad \forall v_h \in V_h. \quad (4.18)$$
$$(u_h, r_H)_\partial - \mu^{-1}(z_H, r_H)_\partial = (\nu^{n+1}, r_H)_\partial, \quad \forall r_H \in Q_0(\partial \mathcal{T}_H), \quad (4.19)$$

where the nudger is given in $z_H = \mu P_H^\partial(u_h - \nu) \in Q_0(\partial \mathcal{T}_H)$ analogous to (4.17) but with $L_H := P_H^\partial$ as defined in (1.16). The nudging scheme (4.18) is subjected first to the saturation test from Section 4.2.1 for the values $\mu = \{1, 4, 64, 16384\}$ of the nudging parameter. Figure 3 demonstrates the resulting evolutions of the error of nudging over time for different test problems. The qualitative discussion of Figure 1 in Section 4.3 can be repeated almost verbatim and the reader is referred to it. We briefly comment on some key differences. As expected from (1.16), the saturation effect manifests with smaller rates (bottom, right). Also, the acceleration of the exponential rate during an initial stage (roughly, for $t \in (0, 0.1)$) is less pronounced. Yet again, in the case of the Kellogg problem (bottom, left) the manifestation of the saturation effect is less clear than for other test problems and the observed rate $\gamma$ (1.11) is substantially higher.

Next we run the convergence test from Section 4.2.2 for values $\ell = \{4, 5, 6, 7\}$ in (4.10). On Figure 4 we present the evolution of errors for test problems from Section 4.1 and the rate of convergence for errors (4.13). The qualitative discussion of Figure 2 in Section 4.3 can be repeated almost verbatim and the reader is referred to it. We briefly comment on some key differences. For the nudging by the boundary projection the observed saturated rate of exponential decay is $\gamma \leq 5$, and the nudging error falls below the discretization level by the time $t = 1$ for the smooth and Dirac problems. Also, the chosen value of nudging parameter $\mu = 64$ is close to the saturated value even in the case of the Kellogg problem $\gamma < 48$. These observations suggest the specific values of $t$ in $L^2(t, T, L^2(\Omega))$ of Theorem 3.10 presented on Figure 4. Note that the boundary projection does not satisfy the stability Assumption 2.2 and, therefore, Proposition 3.5 and Theorem 3.10 do not apply while Proposition 3.2 still guarantees suboptimal rates. Nevertheless, the observed convergence rates are optimal Figure 4.

## 4.5. Nudging by the mean value

Finally, consider the nudging strategy (1.18) with the observation mapping being the mean value, $L_H u := \bar{u} = |\Omega|^{-1} \int_\Omega u$. The full discretization (4.11) of the nudging scheme (1.4) reads: Given $f^{n+1} := f(t_{n+1})$, $g^{n+1} := g(t_{n+1})$, $\bar{\nu}^{n+1} := \bar{\nu}(t_{n+1})$, $n = 0, ..., N-1$, find $u_h := u_h^{n+1} \in V_h$ satisfying for $n = 0, ..., N-1$ the following linear systems,

$$(\partial_\tau u_h, v_h) + a(u_h, v_h) + \mu(\bar{u}_h - \bar{\nu}^{n+1}) \cdot \bar{v}_h = \langle f^{n+1}, v_h \rangle + \langle g^{n+1}, v_h \rangle_\partial, \quad \forall v_h \in V_h. \quad (4.20)$$

The nudging scheme (4.20) is subjected first to the saturation test from Section 4.2.1 for the values $\mu = \{1, 4, 16, 16384\}$ of the nudging parameter. Figure 5 demonstrates the resulting evolutions of the error of nudging over time for different test problems. The qualitative discussion of Figure 1 in Section 4.3 can be repeated almost verbatim and the reader is referred to it. We briefly comment on some key differences. For the smooth (top, left) and the Dirac (top, right) problems the evolution of $\mu = 16384$ has a clear jump during the first time step, i.e. between $u_h(0) = 0$ and $u_h(\tau)$, which can interpreted as follows: the nudging is so effective that the initial stage of nudging is shorter than the step size $\tau$. Interesting enough, the Kellogg problem (bottom, left) is absolutely insensitive to the nudging by the mean value and the evolutions of the nudged scheme correspond to the evolution of the non-nudged scheme (1.8) with the zero initial condition. To explain the situation, we remind that both the exact solution (4.7) and the initial condition $u_h(0) = 0$ have the zero mean value, and the nudging term stays zero for all times.

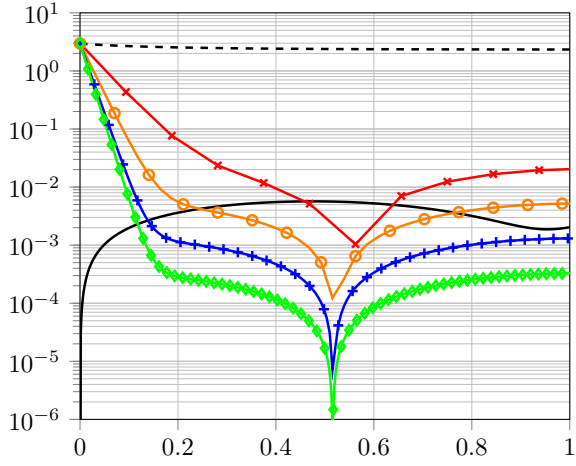

| $\ell$ | $L^2(0.4, 3; L^2)$ | r.o.c |
|---|---|---|
| 4 | 2.327e-2 | - |
| 5 | 5.831e-3 | 2.00 |
| 6 | 1.458e-3 | 2.00 |
| 7 | 3.646e-4 | 2.00 |

| $\ell$ | $H^1(t=3)$ | $L^2(t=3)$ |
|---|---|---|
| 4 | - | - |
| 5 | 1.01 | 1.98 |
| 6 | 1.00 | 2.00 |
| 7 | 1.00 | 2.00 |

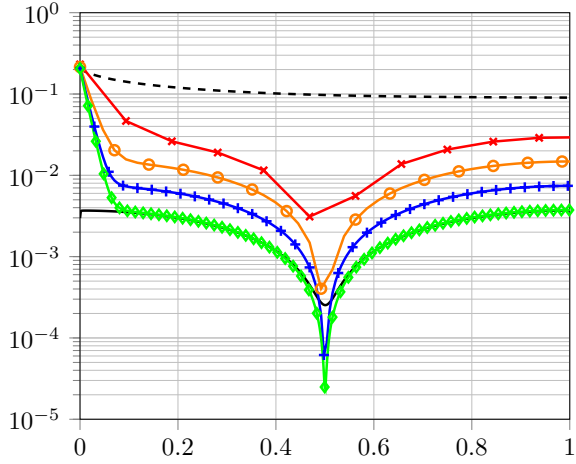

| $\ell$ | $L^2(0.4, 3; L^2)$ | r.o.c |
|---|---|---|
| 4 | 3.364e-2 | - |
| 5 | 1.665e-2 | 1.01 |
| 6 | 8.323e-3 | 1.00 |
| 7 | 4.193e-3 | 0.99 |

| $\ell$ | $H^1(t=3)$ | $L^2(t=3)$ |
|---|---|---|
| 4 | - | - |
| 5 | 0.00 | 1.00 |
| 6 | 0.00 | 0.99 |
| 7 | 0.00 | 0.99 |

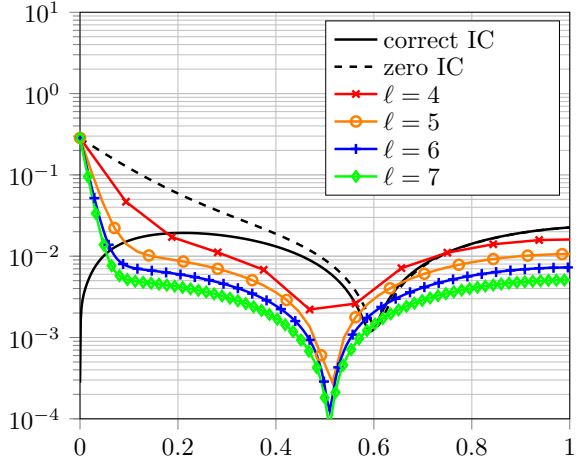

| $\ell$ | $L^2(0.4, 3; L^2)$ | r.o.c |
|---|---|---|
| 4 | 1.841e-2 | - |
| 5 | 1.192e-2 | 0.63 |
| 6 | 8.132e-3 | 0.55 |
| 7 | 5.694e-3 | 0.51 |

| $\ell$ | $H^1(t=3)$ | $L^2(t=3)$ |
|---|---|---|
| 4 | - | - |
| 5 | 0.21 | 0.61 |
| 6 | 0.22 | 0.55 |
| 7 | 0.22 | 0.51 |

Figure 2: Convergence test described in Section 4.2.2 for the nudging by FE projection (1.14). Left: errors (4.12) over time for several $\ell$ in (4.10) for smooth (4.2) (top), Dirac(4.4) (middle) and Kellogg (4.7) (bottom, with the common legend) problems. Right: rates of convergence for errors (4.13). See Section 4.2 for discretization details and Section 4.3 for a discussion of results.

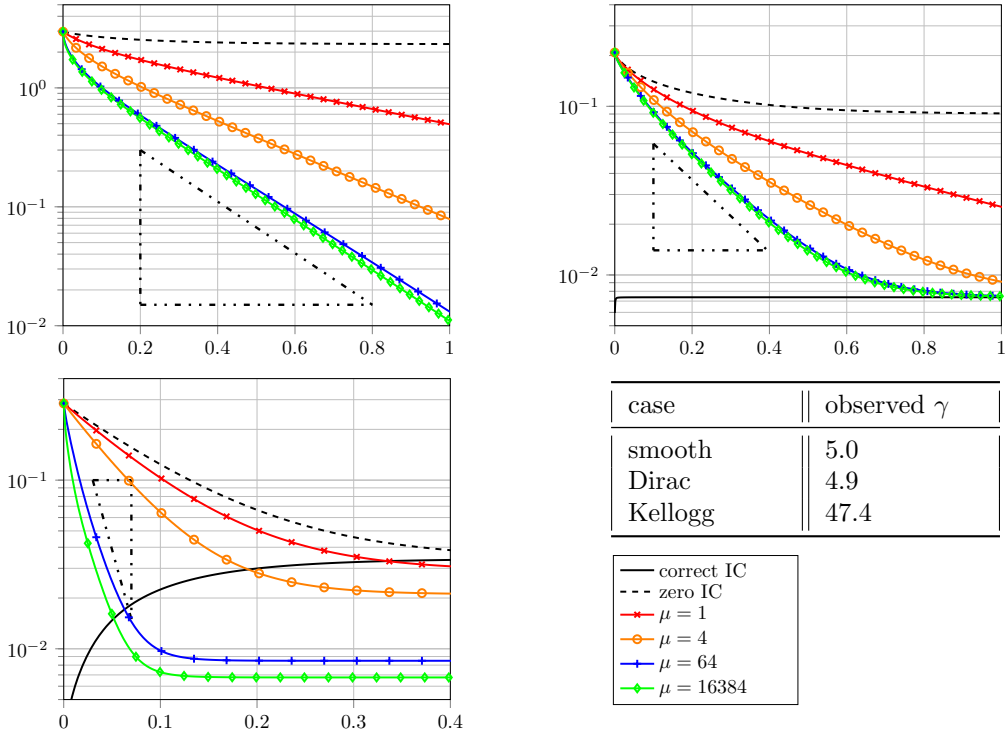

Figure 3: Saturation test described in Section 4.2.1 for nudging by boundary projection (1.16). Errors (4.12) over time for several values of $\mu$ and for the smooth (4.2) (top, left), Dirac (4.4) (top, right) and Kellogg (4.7) (bottom, left) problems. Observed rate $\gamma$ in (1.11) (bottom, right, with the common legend) correspondent to reference triangles. See Section 4.2 for discretization details and Section 4.4 for a discussion of results.

Next we run the convergence test from Section 4.2.2 for $\ell = \{4, 5, 6, 7\}$ in (4.10). On Figure 6 we present the evolution of errors for test problems from Section 4.1 and the rate of convergence for errors (4.13). The qualitative discussion of Figure 2 in Section 4.3 can be repeated almost verbatim and the reader is referred to it. We briefly comment on some key differences. For the nudging by the mean value the observed saturated rates of exponential decay in Figure 5 are $\gamma = 2.3, 2.2, 7.2$, and the nudging error falls below the discretization level by the times $t = 1.7, 2.5, 0.4$ for the smooth, Dirac and Kellogg problems, correspondingly. Also, the chosen value of nudging parameter $\mu = 64$ is greater than the saturated value even in the case of the Kellogg problem $\gamma < 8$. These observations suggest the specific values of $T_\mu$ in $L^2(T_\mu, T, L^2(\Omega))$ of Theorem 3.10 presented on Figure 4. The observed rates of convergence fully support the error analysis of Section 3.

## 5 Conclusion

In this paper an abstract nudging scheme (3.1) suitable for minimally regular problems was studied. The observability property (1.1) granted the exponential decay of the error of nudg-ing on the continuous level as shown in Proposition 2.4. The semi-discrete error analysis of Proposition 3.5 and Theorem 3.10 showed optimal long-term discretization errors for problems with limited or even minimal spatial regularity provided the nudging strategy involved satisfies Assumption 2.1. The numerical experiments in Section 4.2.2 fully supported the optimal error analysis and the saturation effect is demonstrated in Section 4.2.1 in the setting on several nudging schemes and test problems. We conclude the paper with the following possible directions of the future research:

1. Full discretization the nudging schemes (3.1) suitable for problems of limited regularity.

2. The relaxation of Assumption 2.1 for a general nudging strategy as numerical experiments in Section 4.4 demonstrate that it is not necessary for the optimal performance.

3. The analysis of the nudged scheme (3.1) that is augmented by indispensable stabilization forms typical to unfitted FEM and their influence on the long-term error of nudging.

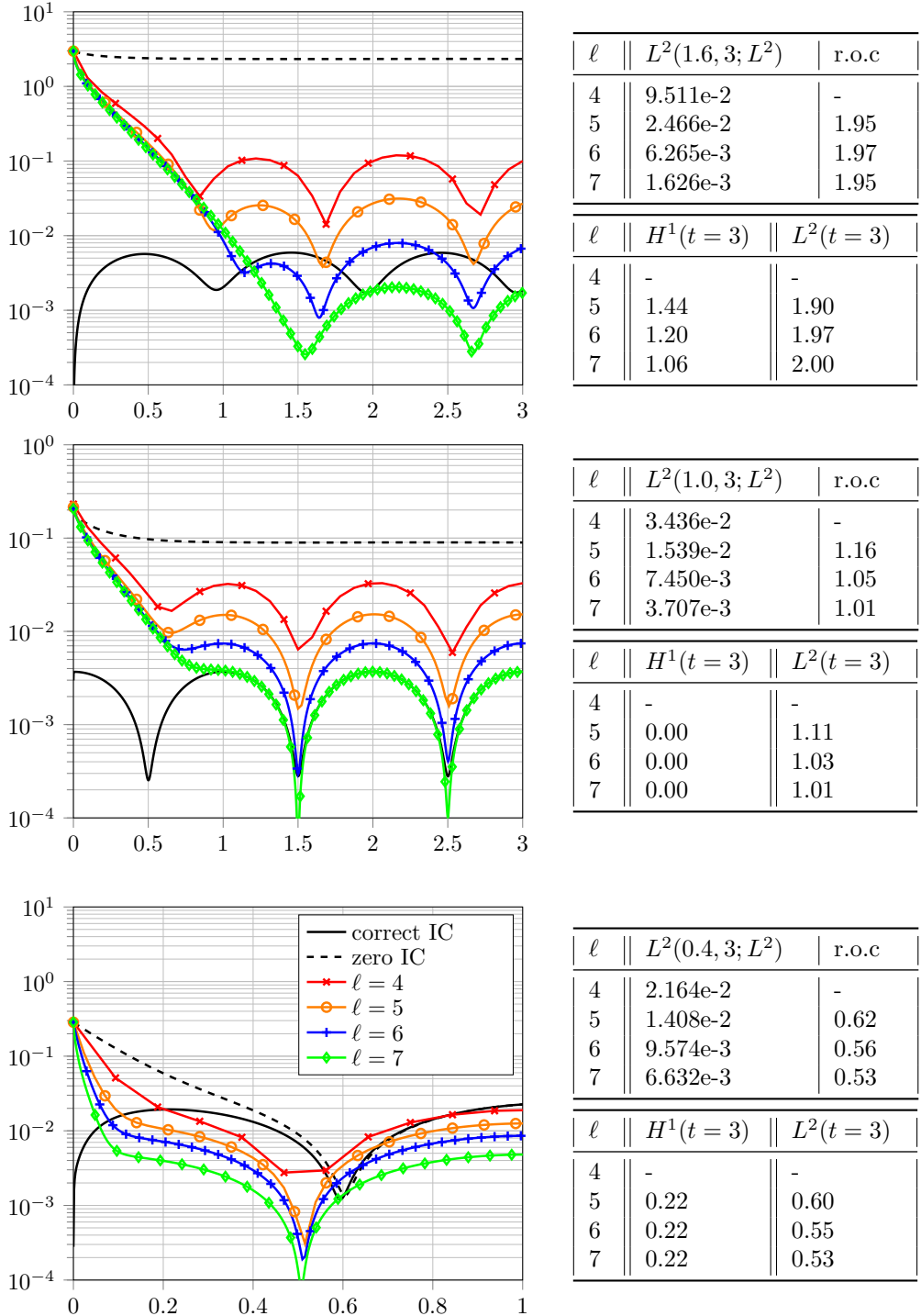

Figure 4: Convergence test described in Section 4.2.2 for nudging by boundary projection (1.16). Left: errors (4.12) over time for several $\ell$ in (4.10) for smooth (4.2) (top), Dirac(4.4) (middle) and Kellogg (4.7) (bottom, with the common legend) problems. Right: rates of convergence for errors (4.13). See Section 4.2 for discretization details and Section 4.4 for a discussion of results.

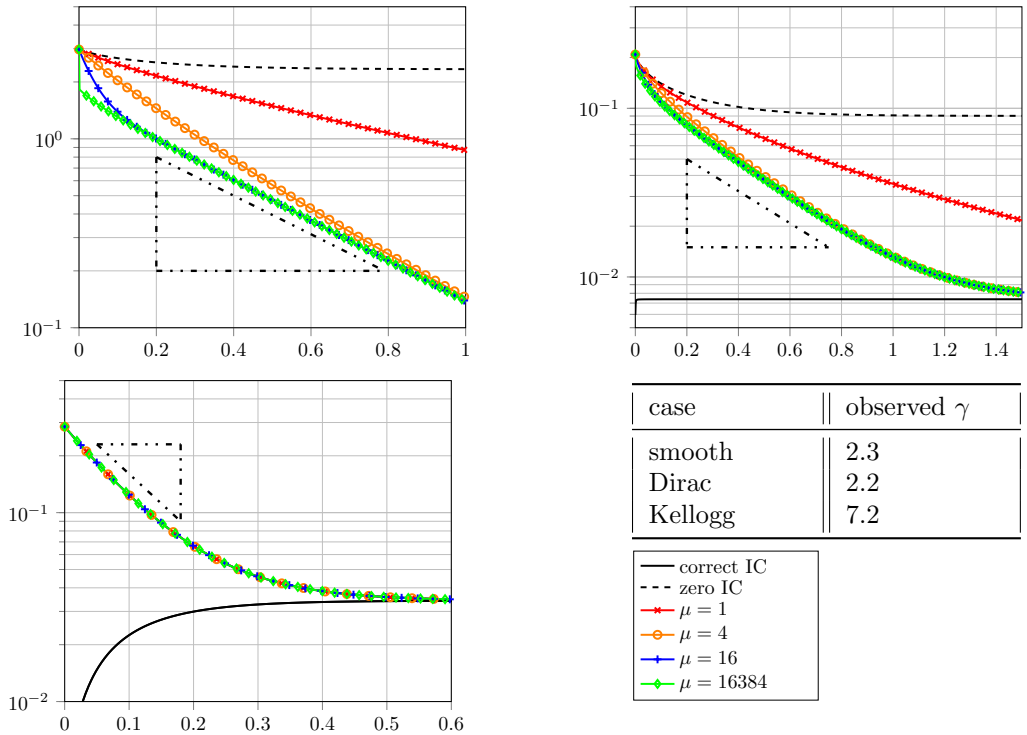

Figure 5: Saturation test described in Section 4.2.1 for the nudging by mean value (1.18). Errors (4.12) over time for several values of $\mu$ and for the smooth (4.2) (top, left), Dirac (4.4) (top, right) and Kellogg (4.7) (bottom, left) problems. Observed rate $\gamma$ in (1.11) (bottom, right, with the common legend) correspondent to reference triangles. See Section 4.2 for discretization details and Section 4.5 for a discussion of results.

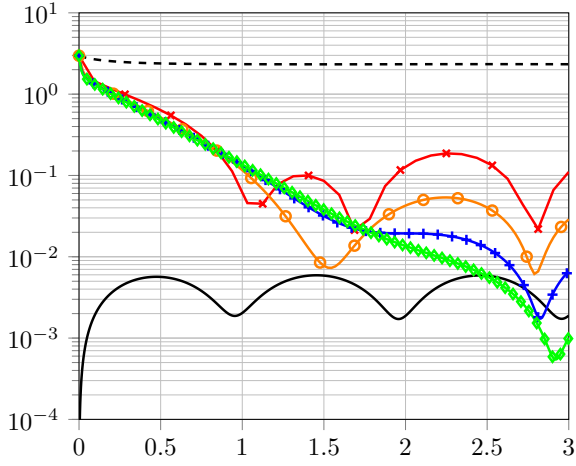

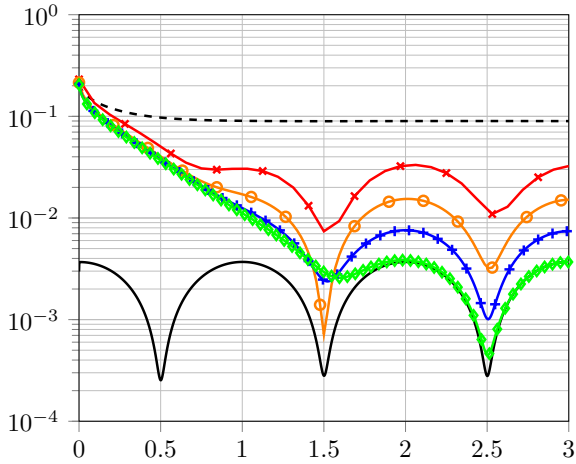

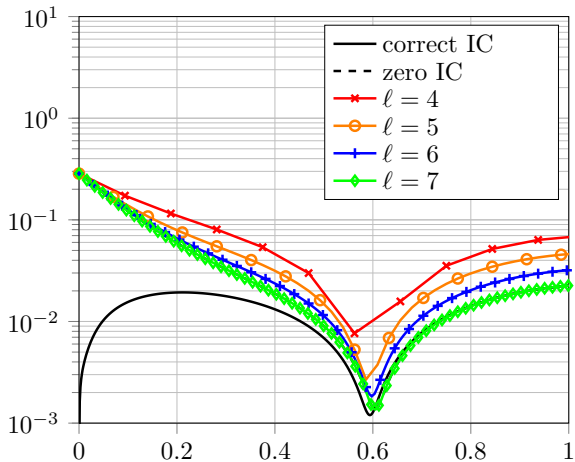

Figure 6: Convergence test described in Section 4.2.2 for nudging by the mean value (1.18). Left: errors (4.12) over time for several $\ell$ in (4.10) for smooth (4.2) (top), Dirac (4.4) (middle) and Kellogg (4.7) (bottom, with the common legend) problems. Right: rates of convergence for errors (4.13). See Section 4.2 for discretization details and Section 4.5 for a discussion of results.